\documentclass[adgeo, manuscript]{copernicus}

\begin{document}
\nolinenumbers

\title{Optimization approaches for the design and operation of open-loop shallow geothermal systems}


\Author[1]{Smajil}{Halilovic}
\Author[2, 3]{Fabian}{Böttcher}
\Author[2]{Kai}{Zosseder}
\Author[1]{Thomas}{Hamacher}

\affil[1]{Chair of Renewable and Sustainable Energy Systems, Technical University of Munich, Garching, Germany}
\affil[2]{Chair of Hydrogeology, Technical University of Munich, Munich, Germany}
\affil[3]{Department for Climate and Environmental Protection (RKU), City of Munich, Germany}




\correspondence{Smajil Halilovic (smajil.halilovic@tum.de)}

\runningtitle{TEXT}

\runningauthor{TEXT}

\received{}
\pubdiscuss{} 
\revised{}
\accepted{}
\published{}


\firstpage{1}

\maketitle

\begin{abstract}
The optimization of open-loop shallow geothermal systems, which includes both design and operational aspects, is an important research area aimed at improving their efficiency and sustainability and the effective management of groundwater as a shallow geothermal resource. 
This paper investigates various approaches to address optimization problems arising from these research and implementation questions about GWHP systems. 
The identified optimization approaches are thoroughly analyzed based on criteria such as computational cost and applicability. 
Moreover, a novel classification scheme is introduced that categorizes the approaches according to the types of groundwater simulation model and the optimization algorithm used.
Simulation models are divided into two types: numerical and simplified (analytical or data-driven) models, while optimization algorithms are divided into gradient-based and derivative-free algorithms.
Finally, a comprehensive review of existing approaches in the literature is provided, highlighting their strengths and limitations and offering recommendations for both the use of existing approaches and the development of new, improved ones in this field.
\end{abstract}

\begin{figure*}[t]
\includegraphics[width=15cm]{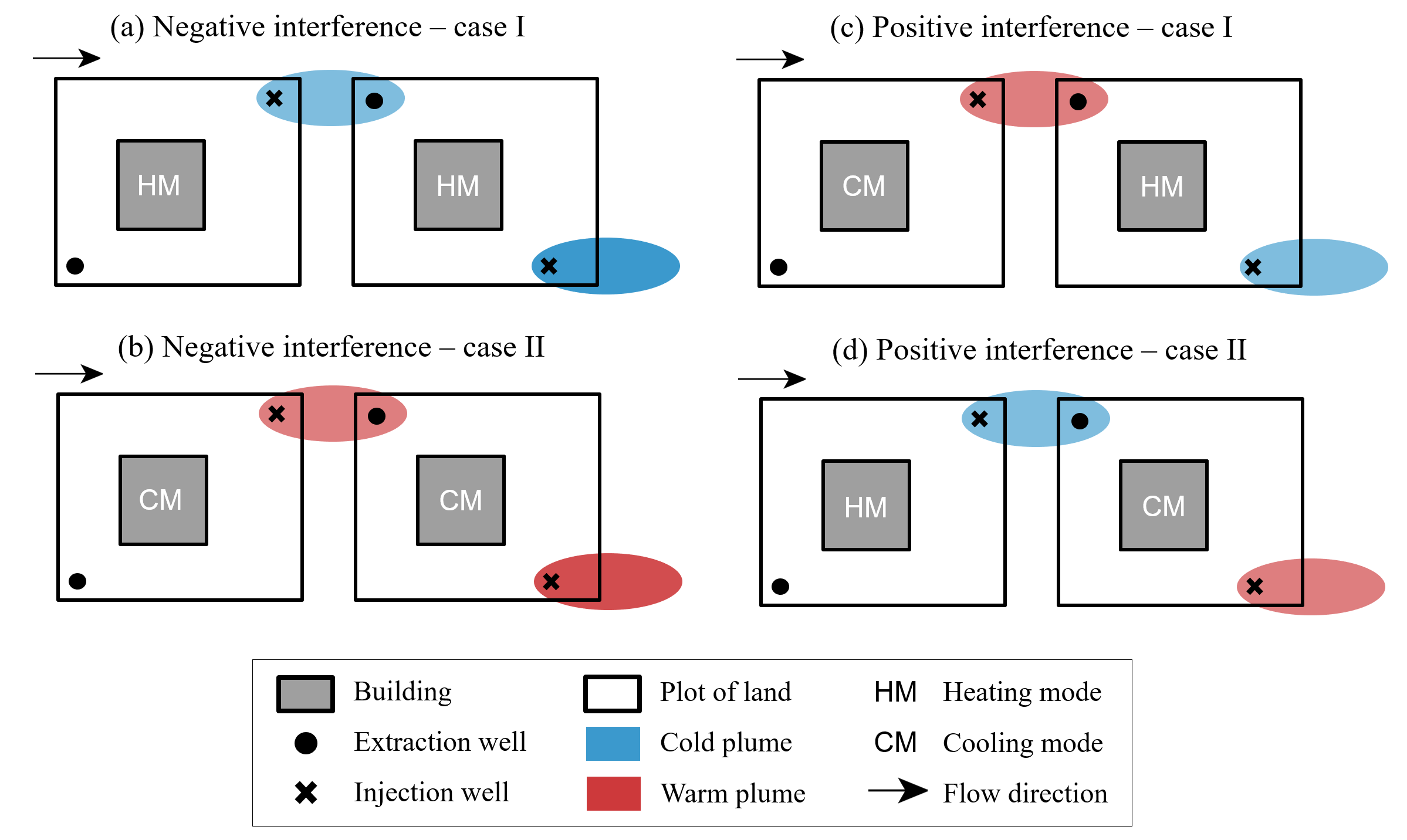}
\caption{Possible thermal interference (interactions) between neighboring GWHP systems.}
\label{fig:interactions}
\end{figure*}

\introduction  
Open-loop shallow geothermal systems, also known as groundwater heat pumps (GWHPs), have emerged as a promising solution for decarbonizing the residential heating and cooling sector \citep{Russo.2012}.
The performance of GWHPs is primarily influenced by groundwater temperature \citep{Kim.2016}, which remains relatively stable throughout the year and is elevated in urban areas due to the subsurface urban heat island effect \citep{Menberg.2013, Epting.2013, Bottcher.2021}.
These systems harness the thermal energy of the aquifer by extracting groundwater from one or more extraction wells and returning it to the same aquifer via injection wells after heat exchange in a heat pump \citep{Florides.2007, Stauffer.2014, GarciaGil.2022}.
Since the temperature of the re-injected water is different from that of the extracted (lower in the heating and higher in the cooling case), this results in thermal plumes in the aquifer that propagate downstream along the groundwater flow direction.
If these plumes reach the extraction wells of neighboring downstream GWHPs, this can result in either negative or positive thermal interference \citep{Perego.2022}.
\autoref{fig:interactions} provides an overview of potential thermal interferences that can occur between neighboring systems, depicting scenarios where the operation of downstream systems can be either degraded (negative interference) or enhanced (positive interference).

It is also important to recognize that GWHPs have a thermal impact on groundwater, which serves as a vital source of drinking water in many places (e.g. \citet{Blum.2021}). 
To mitigate the aforementioned negative interactions and improve the efficiency and sustainability of thermal groundwater use, resource management strategies need to be implemented \citep{Epting.2020}.
This includes optimizing the design, particularly well placement, and operation of GWHP systems, since the propagation of thermal plumes is affected by injection well locations, system operation (pumping rates), and aquifer characteristics.
For example, optimal well placement can minimize negative thermal interference between neighboring GWHPs \citep{Halilovic.2022a} or even maximize their positive interference \citep{GarciaGil.2020}. 
Hence, optimization of multiple neighboring systems plays an important role in urban planning strategies aimed at enhancing sustainability.
In addition, optimization of GWHP systems is crucial for managing groundwater resources and maintaining the current state of groundwater, i.e. preventing adverse changes in its physical, chemical, and biological properties. 
Furthermore, it is essential to ensure adequate spacing between wells within the same GWHP to prevent hydraulic and thermal breakthroughs \citep{Bottcher.2019}. 
Thus, optimization of individual systems is also important to maximize their efficiency and sustainability.
Optimization of individual GWHP systems and concurrent optimization of multiple neighboring systems are challenging due to the complexity of the resulting optimization problems and the necessity for novel and efficient optimization approaches to solve them.

This paper presents a comprehensive overview of optimization approaches for the design and operation of GWHP systems. 
The approaches are critically evaluated and compared based on several criteria, and a novel classification scheme is introduced to effectively categorize these approaches. 
Furthermore, the current status of the approaches found in the literature is presented and possible future research directions are discussed.


\section{Simulation models}\label{sec:simulation}

GWHP systems affect the groundwater body both hydraulically and thermally \citep{GarciaGil.2022}, which can also affect its chemical and biological conditions to a minor degree \citep{Blum.2021}.
The hydraulic head increases around injection wells and decreases around extraction wells, which also changes the hydraulic gradient and groundwater flow patterns.
Thermal impacts are present due to the previously described thermal plumes. 
To analyze these impacts of GWHP systems on groundwater conditions, simulation models are commonly used.
For a particular system design and operation, a simulation model can quantify its impacts on groundwater and, based on that, analyze the performance of the system.
Therefore, simulation models play a crucial role for the computation of the resulting groundwater temperature field and serve as an integral component within the optimization procedures.

These simulation models generally fall into three categories: numerical, analytical, and data-driven.
Numerical models use partial differential equations (PDEs) to describe the underlying physical phenomena, i.e. groundwater flow and heat transport in aquifers.
The resulting system of PDEs can be solved with general PDE solving software or computational fluid dynamics (CFD) software, but there are also several software packages that include specialized domains of numerical simulation for shallow geothermal resources, such as: FEFLOW \citep{Diersch.2014} - based on the finite element method (FEM) or PFLOTRAN \citep{Hammond.2012} - based on the finite volume method (FVM).
Numerical models can incorporate various complex subsurface conditions, including spatially heterogeneous groundwater parameters (e.g. hydraulic conductivity) and conditions (e.g. velocity, temperature, hydraulic head), complex boundary conditions, coupled physical processes, multiple subsurface layers, etc., while simultaneously simulating thermal and hydraulic effects of GWHP systems on the groundwater body.
Therefore, they are closest to reality given sufficient quality of input data, but are generally computationally expensive. 

The second category of models uses analytical formulas to approximate numerical solutions and is commonly applied to estimate thermal plumes associated with smaller GWHPs, whose energy consumption is less than 45,000 kWh per year \citep{Ohmer.2022}.
Due to their analytical nature, these models offer significant computational advantages over numerical models.
In \citet{Pophillat.2020} three prominent analytical models for estimating GWHP thermal plumes were analyzed and compared. 
These models include the radial heat transport model (RHM) \citep{Guimera.2007}, the linear advective heat transport model (LAHM) \citep{Kinzelbach.1987}, and the planar advective heat transport model (PAHM) \citep{Hahnlein.2010}.
The authors concluded that although analytical solutions are less accurate compared to numerical models, they still have value for evaluating the thermal impact of GWHPs. 
Analytical solutions are particularly useful for performing initial assessments of potential negative interference between neighboring GWHPs \citep{Pophillat.2020}.

Finally, data-driven models are gaining popularity in this area of research, primarily due to the emergence of machine learning. 
A common example is the use of neural networks (NNs) to predict thermal plumes \citep{Russo.2014, Leiteritz.2022, Davis.2023}. 
Data-driven models, such as NNs, offer the advantage of fast evaluation, but rely on extensive training data and require additional time for the training process.
Acquiring this training data is often challenging due to the limited measurement and monitoring of hydrogeological data. 
One possible solution is the use of physics-informed neural networks (PINNs) that integrate physical laws driven by PDEs, mitigating the need for extensive training data \citep{Raissi.2017}.


\section{Optimization of GWHPs}\label{sec:optimization}

This section provides a comprehensive analysis of two key aspects related to the optimization of GWHP systems. 
First, in Section 3.1, the underlying optimization problems are discussed. 
Second, in Section 3.2, a detailed overview of the approaches for solving these optimization problems is provided. 
In the following section, we present a generalized problem related to the optimization of GWHP systems, which prepares the way for further analysis in subsequent sections.

\subsection{Optimization problems}\label{ssec:problems}

The high-level optimization problem concerning GWHP systems can be formulated as follows:
\begin{subequations}
\begin{alignat}{2}
&\!\min_{\vec{x}_{\mathrm{d}}, \vec{x}_{\mathrm{o}}} &\qquad& f_{\mathrm{obj}}(\vec{x}_{\mathrm{d}}, \vec{x}_{\mathrm{o}}) = \alpha_1\cdot f_{\mathrm{cost}}(\vec{x}_{\mathrm{d}}, \vec{x}_{\mathrm{o}}) + \alpha_2\cdot f_{\mathrm{env}}(\vec{x}_{\mathrm{d}}, \vec{x}_{\mathrm{o}}) \label{eq:optProb}\\
&\text{subject to} &      & \mathcal{F}_{\mathrm{sim}}(\vec{x}_{\mathrm{d}}, \vec{x}_{\mathrm{o}}) = 0 \label{eq:constraint1}\\
&  &      & \vec{g}(\vec{x}_{\mathrm{d}}, \vec{x}_{\mathrm{o}}) \leq 0 \\
&  &      & \vec{h}(\vec{x}_{\mathrm{d}}, \vec{x}_{\mathrm{o}}) = 0 
\end{alignat}\label{eq:Prob}
\end{subequations} 
\noindent where: 

\begin{tabular}{rl}
$\vec{x}_{\mathrm{d}}$  & vector of optimization variables related to the design of GWHP system(s),\\
$\vec{x}_{\mathrm{o}}$  & vector of optimization variables related to the operation of GWHP system(s),\\
$f_{\mathrm{obj}}$      & objective function to be minimized,\\
$f_{\mathrm{cost}}$     & function describing technical costs,\\
$f_{\mathrm{env}}$      & function describing negative environmental impacts,\\
$\alpha_1, \alpha_2$    & weighting factors,\\
$\mathcal{F}_{\mathrm{sim}}$    & simulation model in a residual form,\\
$\vec{g}$    & inequality constraints,\\
$\vec{h}$    & equality constraints.\\
\end{tabular}
\newline
\newline
\noindent In this generalized problem, we differentiate between two types of optimization variables: design variables $\vec{x}_{\mathrm{d}}$ and operational variables $\vec{x}_{\mathrm{o}}$. 
An example of design variables are the number and spatial layout of GWHP wells, while an example of operational variables are the pumping rates of each well. 
The design variables are constant in time, whereas the operational variables are usually time-dependent.

The objective function $f_{\mathrm{obj}}$ contains two parts: $f_{\mathrm{cost}}$, accounting for the technical costs of GWHP systems, and $f_{\mathrm{env}}$, accounting for the negative environmental impacts. 
The term $f_{\mathrm{cost}}$ can represent various costs associated with the installation and operation of GWHPs, which can be reduced through different means, such as proper sizing of systems (reduced investment costs) or optimal operation of systems (increased efficiency and lifetime). 
On the other hand, the term $f_{\mathrm{env}}$ covers different environmental categories, such as negative impacts on groundwater or $CO_2$ emissions indirectly caused by the operation of GWHP systems. 
It should be noted that environmental considerations are usually incorporated into the problem through constraints and not directly within the objective function. 

The simulation model $\mathcal{F}_{\mathrm{sim}}$ (see \autoref{sec:simulation}) that describes the subsurface phenomena is incorporated into the optimization as a single or multiple equality constraints. 
This model can be of any type discussed previously: numerical, analytical, or data-driven, and it can also have an explicit form, such as PDEs or algebraic equations, or an implicit "black-box" form, such as numerical simulation tools or NNs. 

In addition to the simulation model, other inequality $\vec{g}$ or equality $\vec{h}$ constraints may be present in the optimization problem. 
These can be technical constraints, such as upper and lower limits on pumping rates, regulatory constraints, such as the maximum allowed change in groundwater temperatures, or any other additional constraints.

Depending on how certain elements are specified in the generalized problem \ref{eq:Prob}, the resulting optimization problems can be classified according to different criteria: 
\begin{itemize}
  \item \emph{Optimization variables}: 
  If the only optimization variables are the design parameters $\vec{x}_{\mathrm{d}}$ and the operation of the system(s) is predefined, i.e. $\vec{x}_{\mathrm{o}}$ is fixed, the problem \ref{eq:Prob} becomes a design optimization problem. 
  On the other hand, the problem becomes an optimal control problem when the system design is specified and the operating parameters $\vec{x}_{\mathrm{o}}$ are the optimization variables. 
  Finally, simultaneous optimization of the design and operation of GWHP systems is possible, i.e. considering both $\vec{x}_{\mathrm{d}}$ and $\vec{x}_{\mathrm{o}}$ as optimization variables. 
  This generally leads to improved optimal solutions since there are more degrees of freedom to be optimized, but the resulting problems are usually more difficult to solve due to increased problem complexity (e.g. from linear to non-linear).
  \item \emph{Objective function}: 
  In problem \ref{eq:Prob}, the objective function $f_{\mathrm{obj}}$ is a weighted sum of technical costs $f_{\mathrm{cost}}$ and quantified negative environmental impacts $f_{\mathrm{env}}$. 
  Setting one of the weights $\alpha$ to 0, the problem becomes either a purely economic or an ecological optimization problem, i.e. a single-objective optimization problem.
  If both weighting factors are kept positive, then both the cost and the environmental impact are minimized simultaneously and the problem becomes a type of multi-objective optimization.
  This means that by changing the values of $\alpha_1$ and $\alpha_2$, different Pareto-optimal solutions are obtained \citep{Marler.2010}.

  \item \emph{Application}: 
  The application types can be divided into two main groups: optimization of a single stand-alone system and optimization of multiple neighboring GWHP systems. 
  The application type directly changes the format of the optimization variables and the objective function.
  In addition, different applications may involve different optimization constraints, such as the threshold for negative interference between neighboring systems in the case of optimizing multiple systems. 
  
  \item \emph{Simulation model}:
  In the mathematical sense, the choice of simulation model $\mathcal{F}_{\mathrm{sim}}$ fundamentally changes the type of the optimization problem.
  These optimization problem types belong to different branches of optimization and therefore require corresponding optimization approaches to be solved efficiently.
  For this reason, the entire following section is dedicated to the optimization approaches for the problem \ref{eq:Prob}. 
\end{itemize}

\subsection{Optimization approaches}\label{ssec:approaches}

\begin{figure*}[t]
\includegraphics[width=12cm]{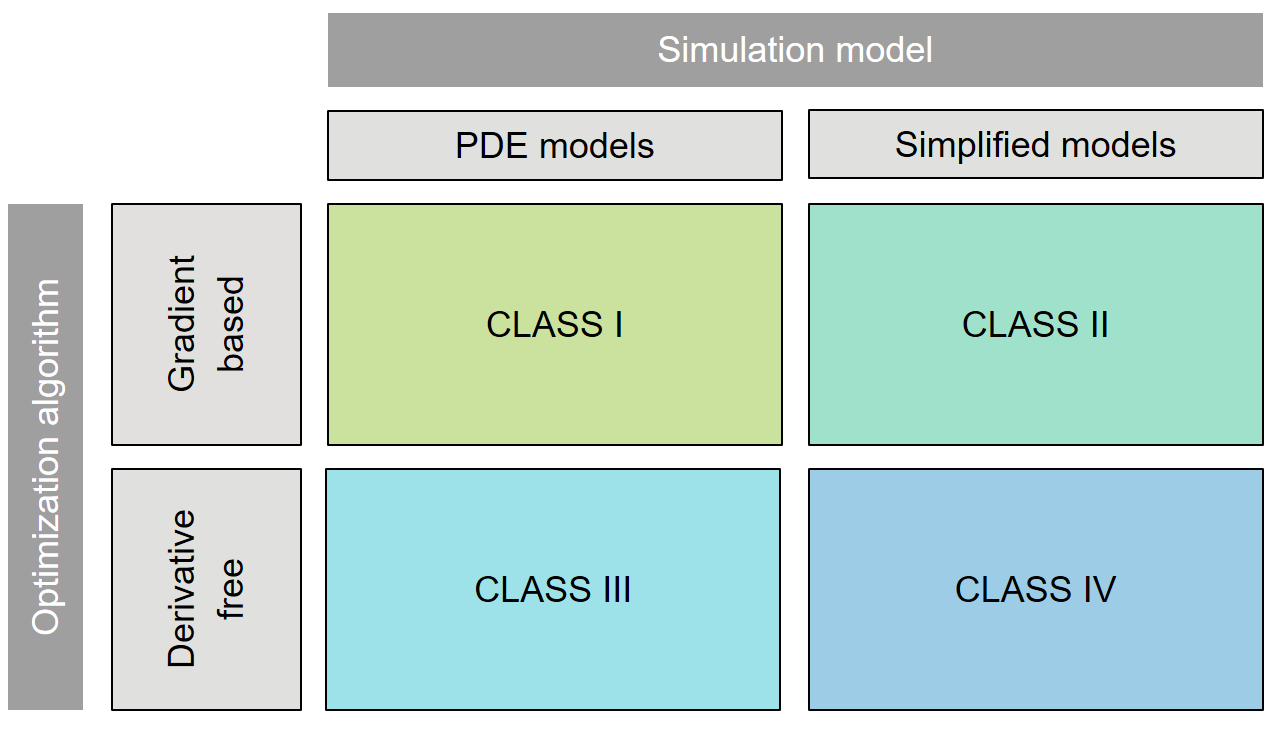}
\caption{Proposed classification of the optimization approaches.}
\label{fig:optim_approaches}
\end{figure*}

In this study, the term "optimization approach" is considered to encompass not only the specific methodology used to solve a given optimization problem, such as the choice of an algorithm, but also the way in which the problem is formulated, which includes the selection of a groundwater simulation model.
The classification of optimization approaches is shown in \autoref{fig:optim_approaches}, where four different classes are identified. 
The categorization is based on the simulation model used, whether it is a PDE model or a simplified model, and the optimization algorithm employed, either gradient-based or derivative-free algorithms. 
In the following, each of these four classes is explained and references to relevant literature sources are provided.

\emph{Class I} comprises optimization approaches where the simulation model is a numerical PDE model, and the optimization is performed using gradient-based algorithms. 
These approaches are referred to as PDE-constrained optimization (PDECO) problems, which are recognized as the most mathematically complex problems of the four classes considered. 
The complexity arises due to the multidisciplinary nature of these problems, necessitating expertise in several areas, including computational optimization, functional analysis, and numerical analysis. 
For example, state-of-the-art groundwater simulation tools usually lack the automatic provision of gradient information, requiring users to estimate gradients manually.
There are two main methods to solve this problem: automatic differentiation of existing simulation tools \citep{Naumann.2011} or the development of custom numerical simulators within frameworks such as Firedrake \citep{Rathgeber.2017} and FEniCS \citep{Logg.2012}, which can automatically provide the required gradients. 
Therefore, a comprehensive understanding of PDE solving is essential in the initial stages of developing a Class I approach.
Various strategies exist for solving PDECO problems, including full-space and reduced-space methods, as well as discretize-then-optimize and optimize-then-discretize approaches \citep{Hinze.2008}. 
For more information on PDECO, the reader is referred to the books by \citet{Hinze.2008} and \citet{Troltzsch.2010}.

\emph{Class II} includes approaches that utilize simplified models for groundwater simulation and employ gradient-based optimization algorithms to solve the underlying optimization problems.
In this context, the simplified models primarily take the form of analytical models (see Section \ref{sec:simulation}). 
These models are expressed through analytical formulas, allowing for direct integration into the optimization problem.
The conceptualization of the overall optimization problem determines the resulting problems, which typically correspond to well-established types of mathematical programming (optimization) problems, such as linear programming (LP), mixed-integer linear programming (MILP), quadratic programming (QP), and similar types. 
These problems are extensively studied in the optimization community, and consequently, efficient algorithms and solvers (software implementations of the algorithms) are readily available.
Comprehensive information on these types of optimization problems can be found in numerous literature sources, including the books by \citet{Nocedal.1999}, \citet{Schrijver.1998} and \citet{Bonnans.2006}.

\emph{Class III} encompasses approaches that combine PDE simulation models with derivative-free optimization (DFO) algorithms.
This form of optimization is commonly referred to as simulation-based optimization, where the simulation model is treated as a black-box, meaning that only the inputs and outputs of the simulator are observed and used by the optimization algorithm to guide the optimization process. 
As a result, the term black-box optimization (BBO) can also be used, which refers to optimization problems where either the objective function or some constraints are treated as black-boxes. 
However, it is important to note that the black-box in BBO is not limited to numerical simulation models. 
For instance, Class IV also falls under the umbrella of BBO. 
Furthermore, the terms BBO and DFO are closely interconnected and can be used interchangeably in certain contexts. 
These distinct terminologies have emerged over time, highlighting different aspects: the conceptual characteristics of BBO, and the algorithmic features of DFO.
DFO algorithms, as the name suggests, do not rely on derivative information during optimization iterations to determine optimal solutions; instead, they use only the values of the objective function and constraints.
These algorithms include heuristic methods (e.g. genetic algorithms, particle swarm optimization, simulated annealing, etc.) \citep{Bozorg.2017}, direct search methods (e.g. MADS algorithm \citep{Le.2011}), and model-based methods (e.g. model-based trust region \citep{Conn.2000}).
For more detailed information on DFO and BBO, interested readers are referred to \citet{Audet.2017} and \citet{Conn.2009}.

\emph{Class IV} comprises optimization approaches that involve the coupling of simplified groundwater simulation models, such as analytical models or NNs, with DFO algorithms. 
Similar to Class III, Class IV falls into the DFO and BBO branches of optimization. Consequently, the same or similar optimization algorithms can be applied to both classes.
The main difference between the two classes lies in the fidelity of the simulation models used. 
Class III uses high-fidelity models, while Class IV relies on low-fidelity models. 
Thus, the computational cost associated with evaluating potential solution candidates during optimization iterations is significantly lower in Class IV.
It is important to note that Class IV has similarities to situations where model-based DFO algorithms are applied in Class III. 
The difference is that the approximate models in Class III are constructed dynamically during the optimization iterations, based on evaluations of the PDE model. 
In contrast, in Class IV, the simplified models are predefined and remain constant throughout the optimization process.

Finally, it should be mentioned that the four introduced classes do not encompass all conceivable approaches, since combined approaches also exist.
For example, the solution obtained from Class II can serve as an initialization for the optimization process in Class I.
However, within the context of this study, the division into four distinct classes seems both logical and practical, since there are substantial differences between these classes.
The following section reviews previous research studies on the optimization of GWHP systems.

\subsubsection{Current status of the approaches used for GWHP systems}\label{ssec:applications}

Despite the increasing importance of GWHP optimization as a research area, the number of existing studies on this topic remains limited.
\citet{Park.2020} propose a simulation-based optimization approach to optimize pumping rates for a single GWHP system. 
The approach couples a numerical groundwater simulation model with a genetic optimization algorithm.
Furthermore, the same approach was extended in \citep{Park.2021} to optimize both well locations and pumping rates within a single system.
Since the approach in these two studies uses a PDE simulation model along with a DFO algorithm, it falls under Class III of the proposed classification.

To date, only one research study has been identified that applies the approach of Class I, i.e. the PDECO framework. 
This study \citep{Halilovic.2022a} introduces a novel approach for concurrently optimizing the well locations of multiple neighboring systems. 
The approach was illustrated using a case study with ten systems, where the optimization objective is to minimize negative interactions between systems and maximize the overall efficiency of all systems.
The proposed approach uses the adjoint method to efficiently compute gradients from the numerical simulation model, which are required by the optimization algorithm. 

There is also only one research study that implements the approach belonging to Class II.
In \citep{Halilovic.2023}, the authors introduce an approach that integrates an analytical groundwater simulation model directly into the optimization problem. 
Specifically, the analytical model used to calculate thermal plumes is the LAHM model \citep{Kinzelbach.1987} and the resulting optimization problem is formulated as an MILP problem.
The study applies this approach to optimize the locations of systems and their associated wells within an urban area comprising 56 potential systems. 
The objective is to satisfy relevant regulations while maximizing heat extraction from the aquifer.
An open-source implementation of the proposed approach can be accessed at \citep{Halilovic.2022d}.

No research studies have been identified that apply the approach of Class IV, which involves the combination of simplified models with DFO algorithms.
It should be noted that other studies on GWHP optimization exist, focusing on aspects such as optimizing the components of a heat pump or determining optimal control strategies. However, these studies do not consider underground processes and are therefore outside the scope of this work.
Furthermore, there are other research studies (e.g. \citet{Zhou.2009, LoRusso.2009, Gao.2013}) that address the optimal design or operation of GWHPs using methods such as scenario comparison or sensitivity analysis. 
However, these methods are not optimization methods, and as such, they are not further discussed in the present study.
In the subsequent section, a comparative analysis is conducted between optimization approaches, i.e. the identified four classes.

\subsubsection{Comparison of the optimization approaches}\label{ssec:comparison}

\begin{figure*}[t]
\includegraphics[width=8.3cm]{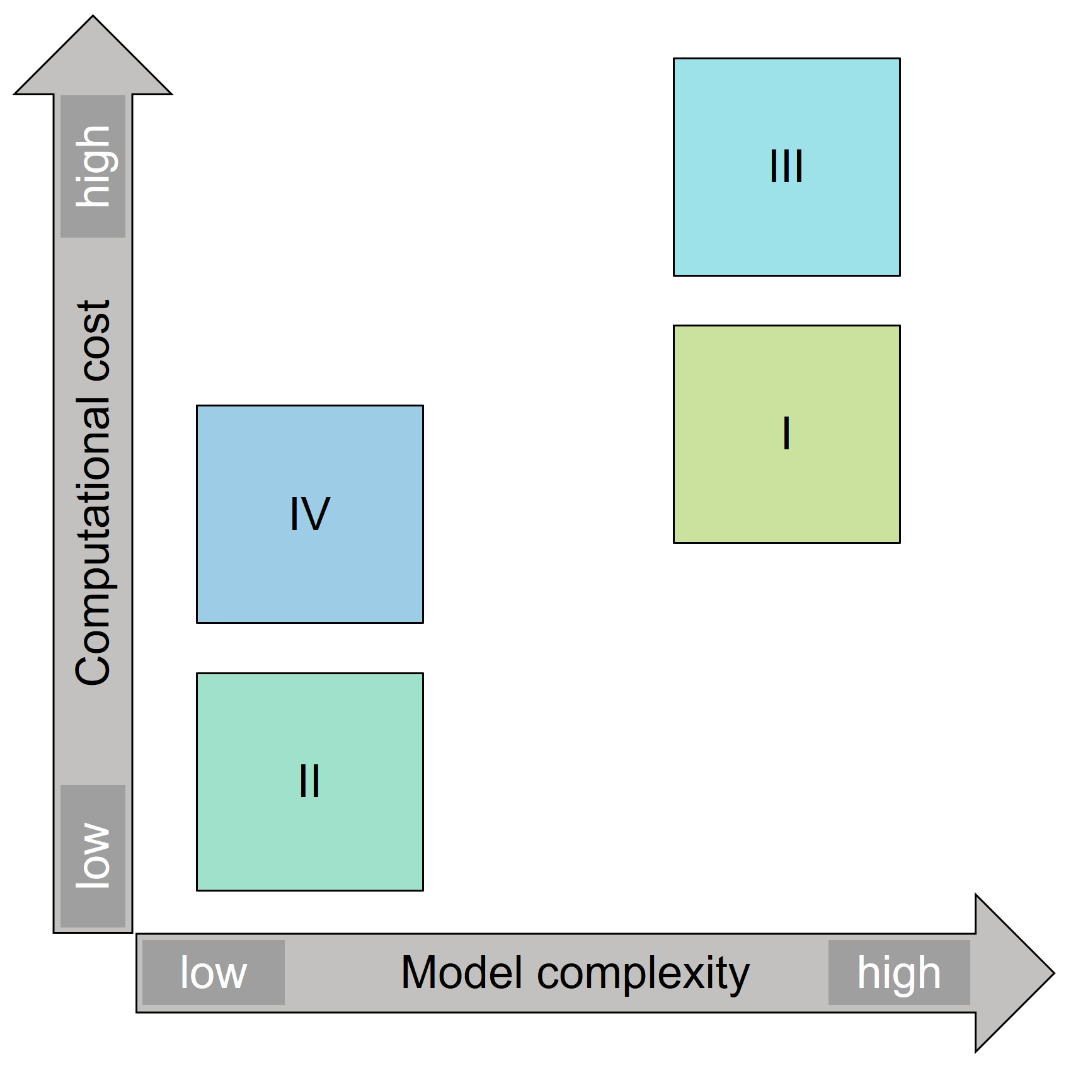}
\caption{Qualitative comparison of the optimization approaches.}
\label{fig:optim_comparison}
\end{figure*}

The primary factors for comparing optimization approaches, i.e. their respective classes, are the computational cost and applicability criteria.
\autoref{fig:optim_comparison} shows a qualitative comparison of these classes, considering two dimensions.
The vertical axis represents the computational cost required to solve optimization problems with the approaches of the respective class.
The computational cost of an approach is of major importance, since in practical planning procedures a relatively fast solution is required.
Moreover, the computational cost increases proportionally with the number of optimization parameters (variables) and the size of the simulation domain. 
Consequently, approaches with high computational costs are limited to scenarios with a small number of optimization parameters and small domains. 
As the number of optimization parameters increases, inefficient approaches quickly become computationally impractical, even when using high-performance computers.

The horizontal axis represents the complexity (fidelity) of the groundwater simulation model used in these approaches.
In the context of this study, the complexity of a simulation model refers to the level of detail in representing physical phenomena in the subsurface, such as the propagation of thermal plumes, that are relevant to the optimization problem under consideration. 
Assuming that the required input data, such as groundwater parameters, are available in sufficient quality, more complex simulation models are more accurate, i.e. closer to reality.
However, it is important to recognize that data on groundwater parameters and conditions are often limited, which limits the use of complex models.
Model complexity is essentially limited by the available data, making the use of highly complex models impractical in the absence of the necessary data. 
Nevertheless, simpler PDE simulation models, such as a 2D model with uniform groundwater conditions, are applicable even with restricted data availability and generally offer higher accuracy than analytical models with identical input data.
Since simulation models are an integral component of optimization approaches, their complexity directly affects the applicability of the obtained optimization results. 
For instance, the results of an approach that uses a complex groundwater simulation model provided with high-quality data can be applied in practice with greater confidence than the results of an approach based on less accurate models.

In the context of computational costs, two key aspects deserve attention: the convergence rate and the computational cost associated with the evaluation of each candidate solution (a unique combination of optimization variables). 
The former quantifies the number of optimization iterations required to reach the optimal solution, while the latter describes the run-time required for each model simulation used to evaluate the current candidate solution within the optimization iterations.
In general, gradient-based algorithms significantly outperform derivative-free algorithms in terms of convergence rate and therefore it is recommended to use gradient-based algorithms when gradient information is readily available and can be obtained at a reasonable cost \citep{Audet.2017, Conn.2009}.
As a result, Class II will almost always outperform Class IV, and Class I will outperform Class III, due to the use of gradient-based algorithms in the former classes (I and II) and derivative-free algorithms in the latter classes (III and IV).
Another disadvantage of Classes III and IV is that derivative-free algorithms generally only find near-optimal solutions and do not guarantee optimality \citep{Audet.2017}.
Furthermore, classes that use simplified simulation models (II and IV) commonly have lower computational costs than classes that use PDE models (I and III).
This is a direct consequence of the computational costs associated with evaluating the simulation model during optimization iterations.
Considering all of the above, a hierarchy of classes based on overall computational cost can be established. 
Class II entails the least computational cost, while Class III is the most computationally demanding.
Classes I and IV fall somewhere in between, with Class IV usually outperforming Class I, although the specific problem characteristics (number of optimization variables, size of the simulation domain, etc.) can also influence this comparative performance.
Despite the computationally demanding nature of Class III, this class is frequently used in the simulation community owing to its user-friendly nature.
By coupling standard standalone simulation software with an existing implementation of a DFO algorithm, typically an evolutionary algorithm, users can develop and apply such approaches relatively quickly.

In terms of the complexity/fidelity of the simulation model used in optimization approaches, it is evident that the classes employing PDE models (I and III) outperform those employing simplified models (II and IV).
The complexity of the simulation model directly influences the validity of the optimization results, thereby affecting the applicability of the corresponding classes. 
Consequently, it is reasonable to use approaches from different classes for different application scenarios. 
For instance, the classes with more complex models (I and III) are suitable for detailed planning of large GWHP systems, while the other two classes (II and IV) can be applied for initial assessments of potential negative interactions between neighboring systems or estimations of geothermal potential on a larger scale.


\section{Discussion and Outlook}\label{sec:discussion}

While there are a limited number of research studies (see Section~\ref{ssec:applications}) that address the optimization of GWHP systems, the research area remains insufficiently explored, which provides an opportunity to pose new research questions and develop novel optimization approaches.
The existing approaches in this field have certain limitations and do not cover all relevant applications.
For example, the approach of Class III presented in \citet{Park.2020, Park.2021} is limited by the number of optimization variables it can efficiently handle, since the computational cost increases exponentially as the number of variables increases.
Similarly, the only study \citep{Halilovic.2022a} using the approach of Class I does not cover all relevant aspects of GWHP optimization, such as optimizing the number of wells in a large GWHP system, optimizing pumping rates, or simultaneously optimizing pumping rates and well locations.

Class I (PDECO) seems to be the most promising among the four classes because it uses PDE simulation models and has lower computational costs compared to Class III. 
Here, the complexity level of the PDE model can be selected based on data availability, as discussed in Section~\ref{ssec:comparison}.
However, this class presents significant challenges due to its mathematical complexity and multidisciplinary nature.
To overcome the challenges associated with developing new approaches within Class I and facilitate their further advancement, collaboration within multidisciplinary teams will be required in the future.

The limitations of the classes that use simplified simulation models (Class II and IV) are directly related to the limitations of the simulation models employed. 
Consequently, improving the simplified simulation models directly enhances the approaches within these classes.
The main goal is to maintain the simulation models as fast and simple to evaluate while enhancing their closeness to reality.
By further improving the accuracy of these simplified models, their scope can be extended to new applications, such as detailed design of large GWHP systems comprising multiple extraction and injection wells.
Moreover, the simplified models are well-suited for integration into energy system optimization models (ESOMs), where GWHP systems play an important role \citep{Halilovic.2022}. 
This is because the coupling of a numerical groundwater simulation with an ESOM is impractical and computationally demanding \citep{Halilovic.2022c}.
By using simplified models, the computational cost can be significantly reduced while still capturing the essential aspects of GWHP systems in the broader context of energy system optimization.
This further enables the development of automated urban planning tools, which will increase sustainability.

Another important consideration in GWHP optimization is the inherent uncertainty associated with subsurface parameters and conditions.
The complex nature of aquifers and the limited availability of measurement and monitoring data contribute to the presence of uncertainties \citep{Gelhar.1992}.
Incorporating these uncertainties into optimization approaches leads to stochastic programming problems \citep{Birge.2011}, which constitute a separate field of optimization. 
The inherent stochastic nature of these problems significantly increases the complexity and computational cost compared to deterministic problems.
Stochastic problems are often solved with modified deterministic optimization approaches or by using a deterministic equivalent of the stochastic problem \citep{Hannah.2015, Li.2021}. 
By minimizing the computational cost of deterministic optimization approaches, researchers can better address the challenges of stochastic problems and develop efficient approaches to such problems. 
Therefore, it is crucial for the deterministic optimization approaches discussed in this study to reduce their computational cost.
Several strategies can be employed to reduce the computational cost of the existing approaches.
For instance, the Class III approach \citep{Park.2020} may improve its efficiency by using model-based algorithms instead of a genetic algorithm. 
Similarly, the Class I approach \citep{Halilovic.2022a} can improve its efficiency by using suitable (re)meshing techniques or by fine-tuning optimization algorithm parameters. 

It is important to note that the classification and comparison of optimization approaches presented in Section~\ref{ssec:approaches} is not limited to the optimization of GWHP systems, but can be applied to any optimization problem where the underlying physical phenomena are described by PDEs.
Moreover, the approaches developed for GWHP systems (see Section~\ref{ssec:applications}) and their future advancements or new approaches in this area can be extended to other applications.
First, they can be extended to applications that share the same underlying physics, such as optimization of aquifer thermal energy storage (ATES) systems, calibration of numerical hydro-thermal groundwater simulation models, and optimization of observation well placement. 
Second, these approaches can be extended to other areas of shallow geothermal energy, including optimization of vertical and horizontal closed-loop shallow geothermal systems and optimization of borehole thermal energy storage (BTES) systems.
Lastly, these approaches can be further extended to areas involving different physical phenomena, such as the optimization of wind farms or tidal power plants.
Nevertheless, it is important to note that advances in these other areas, particularly in the area of shallow geothermal energy, can reciprocally contribute to the improvement of optimization approaches for GWHP systems.
Namely, optimization approaches and principles used in other areas have the potential to be adapted and applied to GWHP systems.

\conclusions\label{sec:conclusions}  

This paper presents a comprehensive analysis and overview of approaches for optimizing the design and operation of GWHP systems. 
First, the optimization problems arising from this research and practice question were investigated, using a generalized problem as a basis.
Then, optimization approaches were identified and compared, and a novel classification of the approaches is proposed. 
The identified approaches were divided into four distinct classes based on the type of groundwater simulation model used (PDE-based or simplified models) and the optimization algorithm applied (gradient-based or derivative-free). 
Finally, the paper includes a thorough review of the existing approaches in the literature, highlighting their limitations and outlining opportunities for future improvements.

Based on the analysis performed, several conclusions can be drawn:
\begin{itemize}
    \item
    Optimization approaches that rely on gradient-based optimization algorithms are preferable, as they consistently outperform derivative-free algorithms.
    \item
    The choice of a simulation model used in an optimization approach has a significant impact on its applicability. For example, approaches using PDE models are more suitable for detailed design of large-scale GWHPs, while simplified models offer practical advantages for assessing the geothermal potential of large areas. However, it is important to note that the degree of model complexity is limited by the availability of hydrogeological data.
    \item
    The existing research on GWHP optimization is limited, with only a few studies addressing this topic.
    \item
    Existing approaches have certain limitations and do not cover all relevant applications and research questions in GWHP optimization.
    One of the main limitations is the high computational cost, which limits the number of optimization parameters and the size of the simulation domain that can be effectively considered. In addition, some approaches are limited in applicability due to the use of simplified groundwater simulation models. Moreover, applications such as optimizing the number and placement of wells in large GWHP systems or simultaneous optimization of pumping rates and well placements remain unexplored.
    Consequently, there is an ongoing need to develop new and improve existing approaches to address these limitations and fill the research gaps.
    \item
    The efficient optimization approaches developed for GWHP systems have the potential to be extended to other shallow geothermal applications as well as to other optimization problems where the underlying physical phenomena are described by PDEs.
    At the same time, approaches from other areas can be adapted and used for GWHP optimization in the future.
\end{itemize}
Finally, this study can provide a valuable foundation for researchers and practitioners involved in the management and optimization of shallow geothermal energy systems.
In particular, it provides valuable insights and recommendations for the application and development of optimization approaches for GWHP systems.
Despite its challenging nature, optimization of GWHP systems is of utmost importance to improve thermal management of groundwater and to unlock the full potential and attractiveness of GWHP technology.

\authorcontribution{SH: Conceptualization, Writing – original draft, Methodology, Investigation, Visualization. FB: Writing - reviewing and editing. KZ: Writing - reviewing and editing. TH: Writing - reviewing and editing, Supervision.} 

\competinginterests{The authors declare no conflicts of interest.} 


\begin{acknowledgements}
No external funding was received for conducting this study.
We would like to thank Jannis Epting for reviewing the paper and providing valuable comments that improved its quality.
\end{acknowledgements}







\bibliographystyle{copernicus}
\bibliography{bibliography.bib}

\end{document}